\documentclass[11pt]{article}
\usepackage{amsfonts,mathrsfs}
\usepackage{amssymb}
\usepackage{amsfonts,amssymb, amsmath, latexsym,theorem}
\newtheorem{theorem}{Theorem}[section]
\newtheorem{proposition}[theorem]{Proposition}
\newtheorem{lemma}[theorem]{Lemma}
\newtheorem{corollary}[theorem]{Corollary}

\newtheorem{remark}[theorem]{Remark}

\def\proof{{\noindent \sc Proof.\quad}}
\newcommand{\proofof}[1]{{\noindent \sc Proof of #1.\quad}}
\def\eproof{{\mbox{}\hfill\qed}\medskip}
\newcommand\qed{{\unskip\nobreak\hfil\penalty50\hskip2em\vadjust{}
  \nobreak\hfil$\Box$\parfillskip=0pt\finalhyphendemerits=0\par}}

\begin{document}

\makeatletter

\renewcommand{\bar}{\overline}

%Other math symbols
\newcommand{\x}{\times}
\newcommand{\<}{\langle}
\renewcommand{\>}{\rangle}
\newcommand{\into}{\hookrightarrow}

%Greek letters
\renewcommand{\a}{\alpha}
\renewcommand{\b}{\beta}
\renewcommand{\d}{\delta}
\newcommand{\D}{\Delta}
\newcommand{\e}{\varepsilon}
\newcommand{\g}{\gamma}
\newcommand{\G}{\Gamma}
\renewcommand{\l}{\lambda}
\renewcommand{\L}{\Lambda}
\newcommand{\n}{\nabla}
\newcommand{\var}{\varphi}
\newcommand{\s}{\sigma}
\newcommand{\Sig}{\Sigma}
\renewcommand{\t}{\theta}
\renewcommand{\O}{\Omega}
\renewcommand{\o}{\omega}
\newcommand{\z}{\zeta}

%Other macros
\newcommand{\p}{\partial}
\renewcommand{\hat}{\widehat}
\renewcommand{\bar}{\overline}
\renewcommand{\tilde}{\widetilde}

%%%%%%%%
% Some fonts
%%%%%%%

%\newcommand{\fiverm}{\tiny\rm}
\font\eightrm=cmr8
\font\ninerm=cmr9

%%%%%
% The next macros define fonts for reals, rationals, complex,
% integers and natural numbers by \R,\Q,\C,\Z and \N respectively.
% Also, \Ri gives \R to the infinity.
%%%%%

\def\N{\mathbb{N}}
\def\Z{\mathbb{Z}}
\def\R{\mathbb{R}}
\def\Q{\mathbb{Q}}
\def\C{\mathbb{C}}
\def\F{\mathbb{F}}
\def\proj{\mathbb{P}}

%%%%%%%
% Some mathematical operators and constants
%%%%%%%

\def\dist{\mathsf{dist}}
\def\Id{{\rm Id}}
\def\Prob{\mathop{\rm Prob}}
\def\Sing{\mathsf{Sing}}
\def\bfE{\mathop{\bf E}}
\def\bd{{\bf d}}
\def\Oh{{\cal O}}
\def\tk{\tilde{\kappa}}
\def\tm{\tilde{\mu}}
\def\dpr{d_{\proj}}
% Para usar en dibujos con Pictex %%
\def\bolita{\scriptscriptstyle\bullet}
\newcommand{\bfx}{{\boldsymbol{x}}}
\def\Disc{{\rm Disc}}
\def\Res{{\rm Res}}

%%%%%%%
% Some abreviations for the bib database
%%%%%%%

\def\JACM{Journal of the ACM}
\def\CACM{Communications of the ACM}
\def\ICALP{International Colloquium on Automata, Languages
            and Programming}
\def\STOC{annual ACM Symp. on the Theory
          of Computing}
\def\FOCS{annual IEEE Symp. on Foundations of Computer Science}
\def\SIAM{SIAM J. Comp.}
\def\SIOPT{SIAM J. Optim.}
\def\BSMF{Bulletin de la Soci\'et\'e Ma\-th\'e\-ma\-tique de France}
\def\CRAS{C. R. Acad. Sci. Paris}
\def\IPL{Information Processing Letters}
\def\TCS{Theoret. Comp. Sci.}
\def\BAMS{Bulletin of the Amer. Math. Soc.}
\def\TAMS{Transactions of the Amer. Math. Soc.}
\def\PAMS{Proceedings of the Amer. Math. Soc.}
\def\JAMS{Journal of the Amer. Math. Soc.}
\def\LNM{Lect. Notes in Math.}
\def\LNCS{Lect. Notes in Comp. Sci.}
\def\JSL{Journal for Symbolic Logic}
\def\JSC{Journal of Symbolic Computation}
\def\JCSS{J. Comput. System Sci.}
\def\JoC{J. Compl.}
\def\MP{Math. Program.}
\sloppy

\bibliographystyle{plain}

%new macros
%%%%%%%%%%%%%%%%%%%

\def\CPRi{{\rm \#P}_{\kern-2pt\R}}

\def\CC{{\mathcal C}}
\def\DD{{\mathcal D}}
\def\NN{{\mathcal N}}
\def\MM{{\mathcal M}}
\def\GG{{\mathcal G}}
\def\UU{{\mathcal U}}
\def\ZZ{{\mathcal Z}}
\def\PP{{\mathscr P}}
\def\AA{{\mathscr A}}
\def\scC{{\mathscr C}}
\def\scE{{\mathscr E}}
\def\sG{{\mathscr G}}
\def\mZ{{\mathcal Z}}
\def\mI{{\mathcal I}}
\def\sH{{\mathscr H}}
\def\sF{{\mathscr F}}
\def\bD{{\mathbf D}}
\def\Hd{\HH_{\mathbf d}}
%%%%%

\def\bx{{\bf x}}
\def\ii{{\'{\i}}}

\def\P{\mathbb P}

\newcommand{\binomial}[2]{\ensuremath{{\left(
\begin{array}{c} #1 \\ #2 \end{array} \right)}}}

\newcommand{\HH}{\ensuremath{\mathcal H}}
\newcommand{\diag}{\mathsf{diag}}
\newcommand{\CH}{\mathsf{CH}}
\newcommand{\Cone}{\mathsf{Cone}}
\newcommand{\SCH}{\mathsf{SCH}}

\newcounter{line}
\newcounter{algorithm}

\newenvironment{algorithm}[3]
{
\addtocounter{algorithm}{1}
{\bf Algorithm \thealgorithm : \sf #1} \\
{\bf Input: }#2 \\
{\bf Output: }#3 \\
\begin{list}{\arabic{line}:}{\usecounter{line}}\setlength{\leftmargin}{3em}
}
{
\end{list}
}
\newcommand{\macheps}{\varepsilon_{\mathrm{m}}}
\newcommand{\sgn}{\mathrm{sgn}}

\begin{title}
{\LARGE {\bf A Numerical Algorithm for Zero Counting. \\
II: Distance to Ill-posedness and Smoothed Analysis}}
\end{title}
\author{Felipe Cucker
\thanks{Partially supported by City University
SRG grant 7002106.}\\
Dept. of Mathematics\\
City University of Hong Kong\\
HONG KONG\\
e-mail: {\tt macucker@cityu.edu.hk}
\and
Teresa Krick
\thanks{Partially supported by grants ANPCyT 33671/05 and UBACyT X113/08-10.}\\
Departamento de Matem\'atica\\
Univ. de Buenos Aires \&\ CONICET\\
ARGENTINA\\
e-mail: {\tt krick@dm.uba.ar}
\and
Gregorio Malajovich\thanks{Partially supported by CNPq grants 
470031/2007-7, 303565/2007-1 and by FAPERJ 
grant E26/110.849/2009.}\\
Depto. de Matem\'atica Aplicada\\
Univ. Federal do Rio de Janeiro\\
BRASIL\\
e-mail: {\tt gregorio@ufrj.br}
\and
Mario Wschebor\\
Centro de Matem\'{a}tica\\
Universidad de la Rep\'{u}blica\\
URUGUAY\\
e-mail: {\tt wschebor@cmat.edu.uy}
}

\date{}
\makeatletter
\maketitle
\makeatother

\begin{quote}
{\small
{\bf Abstract.} We show a Condition Number Theorem for the condition
number of zero counting for real polynomial systems. That is, we show
that this condition number equals the inverse of the normalized distance
to the set of ill-posed systems (i.e.,  those having multiple real zeros).
As a consequence, a smoothed analysis of this condition number follows.
}\end{quote}

\section{Introduction}

This paper continues the work in~\cite{CKMW1}, where we
described a numerical algorithm to
count the number of zeros in $n$-dimensional real projective
space of a system of $n$
real homogeneous polynomials. The algorithm works with
finite precision and both its complexity and the precision
required to ensure correctness are bounded in terms of
$n$, the maximum $\bD$ of the polynomials' degrees,
and a condition number $\kappa(f)$.

In this paper we replace $\kappa(f)$ ---which was originally
defined using the computationally friendly infinity norm---
for a version $\tk(f)$ (defined in Section~\ref{sec:notations}
below) which uses instead
Euclidean norms. This difference is of little consequence
in complexity estimates
since one has (cf.~Proposition~\ref{lem:change} below) that
\begin{equation}\label{eq:ineq}
  \frac{\tk(f)}{\sqrt{n}}\leq\kappa(f)\leq
  \sqrt{2n}\tk(f).
\end{equation}
It allows one, however, to prove a result following
a classical theme in
conditioning ---the relation between condition and
distance to ill-posedness--- and to deduce from this
result a probabilistic analysis for $\tk(f)$ and,
a fortiori, for the complexity of the algorithm
in~\cite{CKMW1}. This analysis is
of interest since, in contrast with
$n$ and $\bD$, the value of $\tk(f)$ (or of $\kappa(f)$)
is not apparent in $f$ and therefore complexity or accuracy
bounds depending on this condition number are not of immediate
use. A solution pioneered by
John von Neumann and collaborators
(see~\cite[\S2.1]{Edelman89} and references therein)
and reintroduced by Steve Smale~\cite{Smale85,Smale97}
is to assume a probability measure on the space of data
and to study the condition number at data $f$
as a random variable. This approach yields bounds
on accuracy or
complexity for random data and has been pursued in several
contexts: systems of linear equations~\cite{Demmel88, Edelman88},
polyhedral conic systems~\cite{ChC01,ChCH05,CW01,HM:06,TTY98},
linear programs~\cite{ChC02}, complex polynomial
systems~\cite{Bez2},~\cite{MALAJOVICH-ROJAS}, etc.
In our case it allows us to trade the presence of $\tk(f)$
in deterministic bounds for probabilistic bounds
in $n$ and $\bD$ only.

\subsection{Distance to Ill-posedness}

It has since long been observed~\cite{Demmel87} that
the condition number
for several problems (in its original aception, as a measure 
of the worst possible
magnification of small input errors in the output~\cite{Rice})
either coincides
with the relativized inverse of the distance from the input to the set of
ill-posed data or is bounded by a small multiple of this inverse. A data 
is {\em ill-posed}
when the magnification above is unbounded. In our case, 
a polynomial system
$f$ is ill-posed when arbitrary small perturbations of $f$ may change its number
of projective real zeros. Systems having this property are exactly those having
multiple (real projective) zeros. Let us denote the set of such systems by
$\Sigma_{\R}$.  Also, let $\Hd$ denote our set of input systems,
i.e. the set of $f=(f_1,\dots,f_n)$, $n$ real homogeneous
polynomials in $n+1$ variables of degrees $\bd:=(d_1,\dots,d_n)$
respectively, endowed with the Bombieri-Weyl norm $\|\ \|_W$ 
(defined in Section~\ref{sec:notations}). Finally, let 
$\dist$ be the distance on $\Hd$ induced by this norm. 

Our main result in this note is the following.

\begin{theorem}\label{main:th}
For all $f\in \Hd$,  $\tk(f)= \frac{\|f\|_W}{\dist(f,\Sigma_{\R})}$.
\end{theorem}

\begin{remark}{\rm
It is worth noting that, although $\tk(f)$ is somehow related to
the condition number $\mu_{\rm norm}(f)$ for complex
polynomial systems (cf.~\cite{BCSS98,Bez1}) a result like
Theorem~\ref{main:th} does not hold for the latter. Actually,
such a result holds on the fibers of the zeros
(see~\cite[Ch.~12, Theorem~3]{BCSS98}) but not globally.

We also note that, as a consequence of Theorem~\ref{main:th},
we have $\tk(f)\geq 1$ for all $f\in\Hd$. This feature, although not
immediate from the definition of $\tk(f)$, follows immediately
from the fact that $0\in\Sigma_\R$.
}
\end{remark}

\subsection{Smoothed Analysis}

Theorem~\ref{main:th} carries out meaningful consequences in
the probabilistic analysis of $\tk(f)$.

The usual probabilistic analysis for polynomial systems
assume that random $f$
are drawn from the unit sphere $S(\Hd)$ with the uniform distribution
(or from such a distribution on the real projective space induced 
by this sphere).
A different approach to the randomization of input data has
been recently
proposed under the name of
{\em smoothed analysis}~\cite{ST:02}. The idea
consists on replacing `random data' by `random perturbations
of given data'. A recent result in~\cite{BuCuLo:07}
derives smoothed analysis bounds for condition numbers
which can be written as inverses  to distances to ill-posedness.
Because of Theorem~\ref{main:th},
these bounds can be straightforwardly applied to $\tk(f)$.

To state this result we need to introduce some notation. Let
$\proj^p(\R)$ denote the  real projective space of dimension $p$ 
and $\dpr$ be the projective distance on $\proj^p(\R)$ (i.e. the sinus
of the Riemannian distance). For a point $a\in\proj^p(\R)$ and
$\sigma\in(0,1]$ we denote by $B(a,\sigma)$ the ball (w.r.t.\
$\dpr$) centered at $a$ and of radius $\sigma$. That is,
$$
   B(a,\sigma):=\{x\in\proj^p(\R)\mid \dpr(x,a)\leq \sigma\}.
$$
In what follows we assume $B(a,\sigma)$ endowed with the uniform
probability measure. Note that in the particular case $\sigma=1$
we obtain $B(a,\sigma)=\proj^p(\R)$ for each $a\in\proj^p$ and 
hence, the usual probabilsitic analysis referred to above.

\begin{theorem}\label{thm:BuCuLo}{\bf \cite{BuCuLo:07}}\quad
Let $S\subset\proj^p$ be contained in
a projective hypersurface $H$ of degree at most $d$ and
$\scC:\proj^p\to[1,\infty]$ be given by
$$
  \scC(z)=\frac{1}{\dpr(z,S)}.
$$
Then, for all $\sigma\in(0,1]$ and all $t\geq (2d+1)\frac{p}{\sigma}$,
\begin{equation*}
  \sup_{a\in\proj^p}\Prob_{z\in B(a,\sigma)}\{\scC(z)\geq t\}\leq
  13dp\frac{1}{\sigma t}
\end{equation*}
and
\begin{equation}\tag*{\qed}
  \sup_{a\in\proj^p}\bfE_{z\in B(a,\sigma)}[\ln \scC(z)]\leq
  2\ln p +2\ln d +\ln\left(\frac{1}{\sigma}\right) +5.
\end{equation}
\end{theorem}

As a consequence of this result we obtain the following corollary
(which we  prove in Section~\ref{sec:proofs}).
Here $\proj(\Hd)$
denotes the projective space associated to $\Hd$,
$N$ denotes its dimension and
$\DD=d_1\cdots d_n$ its associated B\'ezout number.

\begin{corollary}\label{corol:1}
For all $\sigma\in(0,1]$ and all
$t\geq (4n\DD^2+1)\frac{N}{\sigma}$,
\begin{equation*}
  \sup_{f\in\proj(\Hd)}\Prob_{g\in B(f,\sigma)}\{\tk(g)\geq t\}\leq
  13n^2\bD\DD N\frac{1}{\sigma t}
\end{equation*}
and
\begin{equation*}
  \sup_{f\in\proj(\Hd)}\bfE_{g\in B(f,\sigma)}[\ln \tk(g)]\leq
  2\ln N +4\ln n +2\ln \DD +\ln \bD
 +\ln\left(\frac{1}{\sigma}\right) +6.
\end{equation*}
In particular, taking $\sigma=1$, we obtain average analysis:
for all $t\geq N(4n\DD^2+1)$,
\begin{equation*}
  \Prob_{g\in \proj^p}\{\tk(g)\geq t\}\leq
  13n^2\bD\DD N\frac{1}{t}
\end{equation*}
and
\begin{equation*}
  \bfE_{g\in \proj^p}[\ln \tk(g)]\leq
  2\ln N+4\ln n+2\ln \DD+\ln(\bD)+6.
\end{equation*}
\end{corollary}

A recent result~\cite{CHL:08} extends
Theorem~\ref{thm:BuCuLo}
to absolutely continuous measures on $B(a,\sigma)$ whose
densities are radially symmetric around $a$ and may have
a pole at $a$. Applications of this result to $\tk(f)$
readily follow.

\section{Setting and Notations}\label{sec:notations}
For $d\in\N$,  $\HH_d$ denotes the subspace of
$\R[x_0,\ldots,x_n]$ of homogeneous polynomials of degree $d$. We
endow $\HH_d$ with the Bombieri-Weyl inner product $\langle\ ,\
\rangle_W$,  defined for  $f=\sum_{|j|=d}a_j\bfx^j$ and
$g=\sum_{|j|=d}b_j\bfx^j$  by
$$
   \langle f,g\rangle_W=\sum _{|j|=d}\frac{a_jb_j}{{d\choose j}}
$$
where  for $j=(j_0,\dots,j_n)$,  $ {d\choose
j}:=\frac{d!}{j_0!\cdots j_n!}$.
A main feature of this
inner product is its invariance under the action of the orthogonal
group $O(n+1)$. That is, for all $\psi\in O(n+1)$ and all
$f,g\in\HH_d$, $\langle f\circ\psi,g\circ\psi \rangle_W=\langle
f,g\rangle_W$. Next, for $d_1,\ldots,d_n\in\N$, we   endow  $\Hd=\HH_{d_1}\times\ldots\times\HH_{d_n}$ with the inner product 
$$
   \langle f,g\rangle_W=\sum_{i=1}^n\langle f_i,g_i\rangle_W
$$ 
where $f=(f_1,\dots,f_n)$, $g=(g_1,\dots,g_n)\in \Hd$. 
We write $\|\ \|_W$ and $\dist$ to denote 
the norm and distance on $\Hd$ 
induced by this inner product.

Projective zeros of polynomial systems $f \in \Hd$ correspond
to pairs of zeros $(-\zeta,\zeta)$ of the restriction
$f_{|S^n}$ of
$f$ to the $n$-dimensional unit sphere $S^n\subset\R^{n+1}$.
We will thus consider a system $f \in \Hd$ as a (centrally symmetric)
mapping of $S^n$ into $\R^n$.

For a point $x\in S^n$
and a system $f\in\Hd$ one may define both ill-posedness and condition
relative to this point. For the first, one defines
$$
   \Sigma_\R(x)=\{f\in\Hd\mid \mbox{$x$ is a multiple zero of $f$}\},
$$
the set of systems which are ill-posed at $x$. Note that
$\Sigma_\R(x)\neq\emptyset$ for all $x\in S^n$ and that
$$
        \Sigma_\R=\{f\in\Hd\mid \mbox{$f$ has a multiple zero in $S^n$}\}
      =\bigcup_{x\in S^n}\Sigma_\R(x).
$$

Towards the definition of $\tk(f)$,
for $f\in\Hd$ and $x\in S^{n}$, we define
\begin{equation}\label{eq:mu}
  \tm_{\rm norm}(f,x) = \|f\|_{W}\left\| Df(x)_{|T_x S^n}^{-1}
  \left[
  \begin{matrix} \sqrt{d_1} \\ & \sqrt{d_2} \\
   & & \ddots \\ & & & \sqrt{d_n}
  \end{matrix} \right] \right\|
\end{equation}
where $Df(x)_{|T_x S^n}$ is the restriction to the tangent space of
$x$ at $S^n$ of the derivative of $f$ at $x$ and the norm is the
spectral norm, i.e. the operator norm with respect to $\|\ \|_2$.

Next, we define the {\em condition of $f$ relative to $x$} to be
$$
  \tk(f,x)=\frac{\|f\|_W}{\big(\|f\|_W^{2}\tm_{\rm norm}(f,x)^{-2}
  +\|f(x)\|^2_2\big)^{1/2}}.
$$
Finally, we take the {\em condition number} $\tk(f)$ of $f\in \Hd$ to be
its condition relative to its worst conditioned point,
$$
  \tk(f)=\max_{x\in S^n}\tk(f,x)
$$

Note that  for all $\lambda\neq 0$, $\tk(\lambda
f)=\tk(f)$ and $\dist(\lambda
f,\Sigma_\R)=|\lambda|\dist(f,\Sigma_\R)$. The same is true relative to a
point $x\in S^n$. We will therefore assume, without loss of
generality, that $\|f\|_W=1$, and denote by $S(\Hd)$ the unit sphere in
$\Hd$.

\section{The proofs}\label{sec:proofs}

\subsection{The main results}

\begin{proposition}\label{prop:1}
For all $x\in S^n$ and $f\in S(\Hd)$,
$$
   \tk(f,x)=\frac{1}{\dist(f,\Sigma_\R(x))}.
$$
\end{proposition}

\proof For $0\le i\le n$, let $e_i=(0,\ldots,0,1,0,\ldots,0)$ denote
the $i$th coordinate vector. The group $O(n+1)$ acts on $\Hd\times
S^n$ and leaves $\mu_{\rm norm},\tk$ and distance to $\Sigma_\R(\ )$
invariant. Therefore, we may assume without loss of generality that
$x=e_0$. This implies that $T_{e_0}S^n\simeq \langle
e_1,\ldots,e_n\rangle$ and we may write the singular value
decomposition
$$
  \diag\left(\frac{1}{\sqrt{d_i}}\right)Df(e_0)_{|T_{e_0}S^n} =
  \underbrace{\left[\begin{array}{lll}
        & & \\
       u_1 &\ldots & u_n\\
        & &\end{array}\right]}_U
  \left[\begin{array}{lll}
        \sigma_1& & \\
        &\ddots & \\
        & &\sigma_n\end{array}\right]V^{\rm t}
$$
with $U$ and $V$ orthogonal and
$\sigma_1\geq\sigma_2\geq\ldots\geq\sigma_n\geq 0$. Since the
subgroup of $O(n+1)$ leaving $e_0$ invariant is isomorphic to $O(n)$
acting on $T_{e_0}S^n$ we may as well assume that $V=\Id$. Note that
$\tm_{\rm norm}(f,e_0)=\sigma_n^{-1}$, and therefore $ \tk(f,e_0)=
(\sigma_n^2+\|f(e_0)\|_2^2)^{-1/2}$.

Let $g_i(x):=f_i(x)-f_i(e_0)x_0^{d_i}-\sqrt{d_i}\sigma_n
 u_{in}x_0^{d_i-1}x_n$, where $u_n=(u_{1n},\dots,u_{nn})^{\rm t}$.
Clearly, $g_i(e_0)=0$ and $Dg_i(e_0) e_n=0$ (here and in the sequel we
denote $ Dg_i(e_0)_{|T_{e_0}S^n}$ simply by $Dg_i(e_0)$) since
$\partial g_i/\partial x_n (e_0)=\partial f_i/\partial
x_n(e_0)-\sqrt{d_i}u_{in}\sigma_n=0$. Thus,
$g=(g_1,\ldots,g_n)\in\Sigma_{\R}(e_0)$. Moreover
$$
  \|f_i-g_i\|_W^2={d_i\choose d_i}^{-1}f_i(e_0)^2 + {d_i\choose d_i-1,1}^{-1}
  (\sqrt{d_i}\sigma_nu_{in})^2=|f_i(e_0)|^2+\sigma_n^2u_{in}^2
$$
and hence, using $\|u_n\|=1$,
$$
  \|f-g\|_W^2=\|f(e_0)\|_2^2+\sigma_n^2 =\tk(f,e_0)^{-2}.
$$
It follows that
$$
   \dist(f,\Sigma_\R(e_0))\leq \|f-g\|_W =\tk(f,e_0)^{-1}.
$$
\smallskip

For the reciprocal, let $g\in\Sigma_\R(e_0)$. Then, $g(e_0)=0$ and
$Dg(e_0)$ is singular.   We want to show that
$\|f-g\|_W\geq \tk(f,e_0)^{-1}$. To this end, we write
$$
   f_i(x)=f_i(e_0)x_0^{d_i}+\frac{\partial f_i}{\partial x_1}(e_0)x_0^{d_i-1}x_1+\cdots +
   \frac{\partial f_i}{\partial x_n}(e_0)x_0^{d_i-1}x_n+Q_i(x)
$$
with $\deg_{x_0}Q_i\leq d_i-2$ and, similarly,
$$
   g_i(x)=\frac{\partial g_i}{\partial x_1}(e_0)x_0^{d_i-1}x_1+\cdots +
   \frac{\partial g_i}{\partial x_n}(e_0)x_0^{d_i-1}x_n+R_i(x).
$$
Then
$$
  \|f_i-g_i\|_W^2\geq |f_i(e_0)|^2 + \frac{1}{d_i}
         \|Df_i(e_0)-Dg_i(e_0)\|^2_2
$$
and
$$
  \|f-g\|_W^2\geq \|f(e_0)\|_2^2 +
        \left\|\diag\left(\frac{1}{\sqrt{d_i}}\right)Df(e_0)
             -\diag\left(\frac{1}{\sqrt{d_i}}\right)Dg(e_0)\right\|_F^2.
$$
We know that $\diag\left(\frac{1}{\sqrt{d_i}}\right)Dg(e_0)$ is
singular. Hence, denoting by $\Sing_n$ the set of singular
$n\times n$ matrices and by $\dist_F$ the Frobenius distance on
this set of matrices,
$$
  \dist_F\left(\diag\left(\frac{1}{\sqrt{d_i}}\right)Df(e_0),
  \diag\left(\frac{1}{\sqrt{d_i}}\right)Dg(e_0)\right)
  \geq
  \dist_F\left(\diag\left(\frac{1}{\sqrt{d_i}}\right)Df(e_0),
          \Sing_n\right)= \sigma_n,
$$
the last by the Eckart-Young Theorem~\cite[\S11.1]{BCSS98}.
It follows that
\begin{equation}\tag*{\qed}
   \|f-g\|_W^2\geq \|f(e_0)\|_2^2 +\sigma_n^2=\tk(f,e_0)^{-2}.
\end{equation}
\medskip

\proofof{Theorem~\ref{main:th}} Again we can assume $f\in S(\Hd)$.
Note that $$\displaystyle\dist(f,\Sigma_\R)=\min_{g\in\Sigma_\R}\dist(f,g)=%
\min_{x\in S^n}\dist(f,\Sigma_\R(x))$$ since
$\Sigma_\R=\bigcup_{x\in S^n}\Sigma_\R(x)$. Then,
\begin{equation}\tag*{\qed}
   \tk(f)=\max_{x\in S^n}\tk(f,x)=
    \max_{x\in S^n}\frac{1}{\dist(f,\Sigma_\R(x))}
    =\frac{1}{\displaystyle\min_{x\in S^n}\dist(f,\Sigma_\R(x))}
    =\frac{1}{\dist(f,\Sigma_\R)}.
\end{equation}

Before proving Corollary~\ref{corol:1} we recall some useful facts
in algebraic geometry.

For $1\le i\le n$, let $f_i= \sum_{|j|=d_i}u_{ij}\bfx^j$ be a
generic (i.e. with indeterminate coefficients) homogeneous
polynomial of degree $d_i$ in the variables $\bfx=(x_0,\dots,x_n)$
and $f=(f_1,\dots,f_n)$. Set $N:=\sum_{i=1}^n {d_i+n\choose n} -1$,
the dimension of the projective  coefficients space. The
$\mathbf{d}$-{\em discriminant variety} $\Sigma_\C\subset \P^N(\C)$
is the locus of such polynomial systems $f=(f_1,\dots,f_n)$ with
multiple zeros, i.e. such that there exists $z \in\C^{n+1}$,
$z\neq 0$, with $f_1(z) = \cdots = f_n(z) = 0
$ and $Df(z)$ has rank $< n$. It is well-known that $\Sigma_\C$ is a
hypersurface in $\P^N(\C)$ defined by an irreducible polynomial
$\Disc(f)\in \Z[u_{ij}]$ (see \cite{Jou} or \cite[Ch.~10]{GKZ}). For
lack of a precise reference we prove the following result.

\begin{lemma}\label{lem:degree}
 $$\deg(\Sigma_\C)=n{\cal D} + (d_1+\cdots +d_n-n-1){\cal D}
\sum_{j=1}^n\frac{1}{d_j}.$$
\end{lemma}

\proof We know that $\deg(\Sigma_\C)=\deg(\Disc(f))$. We apply
Identity~(13) of \cite{DanJer}:
$$
     \Res_{\rho,d_1,\dots,d_n}(J_f,f_1,\dots,f_n)
   = \Res_{d_1,\dots,d_n}(f_1^0,\dots,f_n^0)\Disc(f),
$$
where the standard notation $\Res_{d_1,\dots,d_n}$ is for the
multihomogeneous projective resultant of $n$ generic homogeneous
polynomials of respective degrees $d_1,\dots,d_n$ in $n$ variables,
$\rho:= d_1+\cdots +d_n-n$, $J_f$ is the determinant of the matrix
$(\partial f_i/\partial x_j)_{1\le i,j\le n}$ and $f_i^0$ denotes
the homogeneous component (of degree $d_i$) of
$f_i(1,x_1,\dots,x_n)$.

We note that $\deg(J_f)$ is a polynomial of degree $\rho$ in $\bfx$
whose coefficients are polynomials in  $u_{ij}$ of degree $n$.
 On the
other hand $\Res_{d_1,\dots,d_n}$ is a multihomogeneous polynomial
of degree $\prod_{k\ne i} d_k$ in the group of variables $u_{ij}$
(\cite{CLO}). Therefore, since $J_f$ has degree $n$ in the $u_{ij}$,
we derive
$$
    n{\cal D} + \rho \sum_{1\le j\le n}\frac{\cal D}{d_j}
   =\sum_{1\le j\le n}\frac{\cal D}{d_j} + \deg\big(\Disc(f)\big).
$$
The statement easily follows.
 \eproof

\proofof{Corollary~\ref{corol:1}} We have  $\Sigma_{\R}\subset
\Sigma_\C$ and for all $f\in\Hd$, $f\neq 0$,
$$
   \tk(f) = \frac{\|f\|_W}{\dist(f,\Sigma_{\R})}
              =\frac{1}{\dpr(f,\Sigma_{\R}\cap \P^N(\R))}.
$$
The statement now follows from Theorem~\ref{thm:BuCuLo} with $p=N$,
$S=\Sigma_{\R}\cap \P^N(\R)$, $H=\Sigma_{\C}\cap \P^N(\R)$,
and $d=n^2\bD\DD$, the last
by~Lemma~\ref{lem:degree} since $n{\cal D} + (\rho-1){\cal D}
\sum_{j=1}^n\frac{1}{d_j} \le n^2\bD\cal D$. \eproof

\subsection{The equivalence of condition numbers}

In~\cite{CKMW1} we considered $\Hd$ endowed with the norm
given by
$$
    \|f\|:=\displaystyle{\max_{1\le i\le n}\|f_i\|_W}
$$
(we note $\|f\|\le \|f\|_W\le \sqrt n \|f\|$) and defined
$$
   \kappa(f):=\max_{x\in S^n} \min
     \left\{\mu_{\rm norm}(f,x),\frac{\|f\|}{\|f(x)\|_\infty}\right\}
$$
with
$$
   \mu_{\rm norm}(f,x)
  :=\sqrt{n}\,\|f\|\, \left\| Df(x)_{|T_x S^n}^{-1}
  \diag (\sqrt{d_i})
 \right\|= \sqrt n \frac{\|f\|}{\|f\|_W} \tm_{\rm norm}(f,x)
$$
(we note $\tm_{\rm norm}(f,x)\le \mu_{\rm norm}(f,x)\le\sqrt n \,\tm_{\rm norm}(f,x)$).
The next result shows that the $\kappa(f)$ thus defined is
closely related to $\tk(f)$.

\begin{proposition}\label{lem:change}
$$
  \frac{\tilde\kappa(f)}{\sqrt n} \le \kappa(f)\le
   \sqrt{2n}\ \tilde\kappa(f).
$$
\end{proposition}

\proof Let $x\in S^n$.
We observe that since $\|f\|_W^{2}\tm_{\rm norm}(f,x)^{-2}+ \|f(x)\|_2^2\ge \|f\|_W^{2}\tm_{\rm norm}(f,x)^{-2}$,
we have
$\tk(f,x) \le \tm_{\rm norm}(f,x)\le \mu_{\rm norm}(f,x)$. Similarly, using $\|f\|_W\le \sqrt n \|f\|$ and $\|f(x)\|_2\ge \|f(x)\|_\infty$,
we have $\tk(f,x) \le \sqrt n \|f\|/\|f(x)\|_\infty$. Therefore
$$
    \tk(f,x)\le \sqrt n\min
  \left\{ \mu_{\rm norm}(f,x),\frac{\|f\|}{\|f(x)\|_\infty}\right\},
$$
which implies
$$
   \tk(f)=\max_{x\in S^n}\tk(f,x)\le 
  \sqrt n \max_{x\in S^n}\min
  \left\{ \mu_{\rm norm}(f,x),\frac{\|f\|}{\|f(x)\|_\infty}\right\}
 =\sqrt n \,\kappa(f).
$$

To prove the other inequality note that, for any $x\in S^n$,
\begin{eqnarray*}
  \min\left\{\|f\|^2\left\| Df(x)_{|T_x S^n}^{-1}
  \diag (\sqrt{d_i})
 \right\|^{2},
  \frac{\|f\|^2}{\|f(x)\|_2^2}\right\}&\le &
  \frac{2\|f\|^2}{\left\| Df(x)_{|T_x S^n}^{-1}
  \diag (\sqrt{d_i})
 \right\|^{-2}+\|f(x)\|_2^2}\\
 &\le&  \frac{2\|f\|_W^2}{\|f\|_W^2 \tm_{\rm norm}(f,x)^{-2}
 +\|f(x)\|_2^2}\\[2mm]
 &=& 2\tk(f,x)^2,
\end{eqnarray*}
and $\|f(x)\|_2\le \sqrt n \, \|f(x)\|_\infty$. Therefore,
$$
   \min\left\{ \mu_{\rm norm}(f,x), \frac{\|f\|}{\|f(x)\|_\infty}\right\}
   \le \sqrt{2n}  \,\tk(f,x).
$$
This  implies
$
  \kappa(f)\le 
 \sqrt{2n} \, \tk(f)
$.
\eproof

\noindent{\bf Acknowledgment.}\quad
In June 2008 we discussed condition and distance to ill-posedness
with Mike Shub who suggested to us that a result in the spirit of
Theorem~\ref{main:th} should be true. We are grateful to him
since that conversation is at the origin of our paper.

{\small
%   \bibliography{../../../book/book}

\begin{thebibliography}{1}

\bibitem{BCSS98}
L.~Blum, F.~Cucker, M.~Shub, and S.~Smale.
\newblock {\em Complexity and Real Computation}.
\newblock Springer-Verlag, 1998.

\bibitem{BuCuLo:07}
P.~B\"urgisser, F.~Cucker, and M.~Lotz.
\newblock The probability that a slightly perturbed numerical analysis problem
  is difficult.
\newblock {\em Mathematics of Computation}, 77:1559--1583, 2008.

\bibitem{ChC01}
D.~Cheung and F.~Cucker.
\newblock Probabilistic analysis of condition numbers for linear programming.
\newblock {\em Journal of Optimization Theory and Applications}, 114:55--67,
  2002.

\bibitem{ChC02}
D.~Cheung and F.~Cucker.
\newblock Solving linear programs with finite precision: {I}. {C}ondition
  numbers and random programs.
\newblock {\em \MP}, 99:175--196, 2004.

\bibitem{ChCH05}
D.~Cheung, F.~Cucker, and R.~Hauser.
\newblock Tail decay and moment estimates of a condition number for random
  linear conic systems.
\newblock {\em \SIOPT}, 15:1237--1261, 2005.

\bibitem{CLO}
D.~Cox, J.~Little, and D.O'Shea. \newblock {\em Using algebraic
geometry.}
\newblock Graduate Texts in
Mathematics, 185. Springer-Verlag, New York, 1998.

\bibitem{CHL:08}
F.~Cucker, R.~Hauser, and M.~Lotz.
\newblock Adversarial smoothed analysis.
\newblock Preprint, 2008.

\bibitem{CKMW1}
F.~Cucker, T.~Krick, G.~Malajovich, and M.~Wschebor.
\newblock A numerical algorithm for zero counting. {I}: {C}omplexity and
  accuracy.
\newblock {\em \JoC}, 24:582--605, 2008.

\bibitem{CW01}
F.~Cucker and M.~Wschebor.
\newblock On the expected condition number of linear programming problems.
\newblock {\em Numer. Math.}, 94:419--478, 2003.

\bibitem{DanJer}
C.~D'Andrea, and  G.~Jeronimo. \newblock Rational Formulas for
Traces in zero-dimensional Algebras.
\newblock Preprint, 2008.

\bibitem{Demmel87}
J.~Demmel.
\newblock On condition numbers and the distance to the nearest ill-posed
  problem.
\newblock {\em Numer. Math.}, 51:251--289, 1987.

\bibitem{Demmel88}
J.~Demmel.
\newblock The probability that a numerical analysis problem is difficult.
\newblock {\em Math. Comp.}, 50:449--480, 1988.

\bibitem{Edelman88}
A.~Edelman.
\newblock Eigenvalues and condition numbers of random matrices.
\newblock {\em SIAM J. of Matrix Anal. and Applic.}, 9:543--556, 1988.

\bibitem{Edelman89}
A.~Edelman.
\newblock {\em Eigenvalues and Condition Numbers of Random Matrices}.
\newblock Ph. D. Thesis, M.I.T., 1989.

\bibitem{HM:06}
R.~Hauser and T.~M\"uller.
\newblock Conditioning of random conic systems under a general family of input
  distributions.
\newblock {\em Found. Comput. Math.}, 9:335--358, 2009.

\bibitem{GKZ}  I.M.~Gel'fand, M.M.~Kapranov, and A.V.~Zelevinsky.
\newblock {\em Discriminants, resultants, and multidimensional
determinants.}
\newblock
Mathematics: Theory \& Applications. Birkh¨auser Boston, Inc.,
Boston, MA, 1994. x+523 pp.

\bibitem{Jou} J.P.~Jouanolou. \newblock Cours DEA,
\newblock University of Strasbourg.

\bibitem{MALAJOVICH-ROJAS}
G.~Malajovich and J.~M. Rojas.
\newblock High probability analysis of the condition number of sparse
  polynomial systems.
\newblock {\em \TCS}, 315:524--555, 2004.

%\bibitem{MORGAN-SOMMESE}
%A. Morgan and A. Sommese.
%\newblock  A homotopy for solving general polynomial
%systems that respects $m$-homogeneous structures.
%\newblock  {\em Appl. Math. Comput.}, 24:101--113, 1987.

\bibitem{Rice}
J.R. Rice.
\newblock A theory of condition.
\newblock {\em SIAM J. Numer. Anal.}, 3:217--232, 1966.

%\bibitem{SHAFAREVICH}
%I. R. Shafarevich.
%\newblock {\em Basic algebraic geometry}.
%Springer Study Edition, Springer-Verlag, Berlin, 1977.

\bibitem{Bez1}
M.~Shub and S.~Smale.
\newblock Complexity of {B}\'ezout's theorem {I}: geometric aspects.
\newblock {\em \JAMS}, 6:459--501, 1993.

\bibitem{Bez2}
M.~Shub and S.~Smale.
\newblock Complexity of {B}\'ezout's theorem {II}: volumes and probabilities.
\newblock In F.~Eyssette and A.~Galligo, editors, {\em Computational Algebraic
  Geometry}, volume 109 of {\em Progress in Mathematics}, pages 267--285.
  Birkh\"auser, 1993.

\bibitem{ST:02}
D.A. Spielman and S.-H. Teng.
\newblock Smoothed analysis of algorithms.
\newblock In {\em Proceedings of the International Congress of Mathematicians},
  volume~I, pages 597--606, 2002.

\bibitem{Smale85}
S.~Smale.
\newblock On the efficiency of algorithms of analysis.
\newblock {\em \BAMS}, 13:87--121, 1985.

\bibitem{Smale97}
S.~Smale.
\newblock Complexity theory and numerical analysis.
\newblock In A.~Iserles, editor, {\em Acta Numerica}, pages 523--551. Cambridge
  University Press, 1997.

\bibitem{TTY98}
M.J. Todd, L.~Tun\c{c}el, and Y.~Ye.
\newblock Characterizations, bounds and probabilistic analysis of two
  complexity measures for linear programming problems.
\newblock {\em \MP}, 90:59--69, 2001.

\end{thebibliography}

}

\end{document}